\documentclass[11 pt]{amsart}
\usepackage{amsmath,enumerate,amsthm,amssymb,amscd}

\newcommand{\ra}{\rightarrow}

\newcommand{\gothic}{\mathfrak}

\newcommand{\m}{{\gothic{m}}}

\newcommand{\depth}{\operatorname{depth}}

\renewcommand{\hat}{\widehat}
\renewcommand{\phi}{\varphi}

\newcommand{\Ass}{\operatorname{Ass}}

\newcommand{\grade}{\operatorname{grade}}

\newcommand{\lR}{\lambda_{R}}

\newtheorem{thm}{Theorem}
\newtheorem{cor}[thm]{Corollary}
\newtheorem{prop}[thm]{Proposition}
\newtheorem{lemma}[thm]{Lemma}
\newtheorem{defn}[thm]{Definition}
\newtheorem{remark}[thm]{Remark}
\newtheorem{ques}[thm]{Question}
\newtheorem{example}[thm]{Example}

\numberwithin{equation}{section}

\begin{document}
\title{Hilbert Coefficients of Parameter Ideals}

\author {Lori McCune}

\address{Department of
Mathematics and Computer Science\\
Ashland University\\Ashland,  OH 44805}
\email{lmcdonne@ashland.edu}

\date{\today}
\
\keywords{Hilbert-Samuel polynomial, Hilbert coefficients}
\subjclass[2000] {Primary 13D40}
\bibliographystyle{amsplain}

\numberwithin{thm}{section}

\begin{abstract} We consider the non-positivity of the Hilbert coefficients for a parameter ideal of a commutative Noetherian local ring.  In particular, we show that the second Hilbert coefficient of a parameter ideal of depth at least $d-1$ is always non-positive and give a condition for the coefficient to be zero.  With the added condition that the depth of the associated graded ring is also at least $d-1$ we show $e_{i}(q) \leq 0$ for $i=1, \dots, d$.
\end{abstract}

\maketitle

\section{Introduction}

Let $(R, \m)$ be a Noetherian local ring, $I \subseteq R$ an $\m$-primary ideal and $M$ a finitely generated $R$-module.  We let $\lR(-)$ denote the length of an $R$-module. The \emph{Hilbert function} for $I$ with respect to $M$ is the function $H_{I,M}: \mathbb{Z} \rightarrow \mathbb{Z}$ given by $H_{I,M}(n)=\lR(M/I^{n}M)$.   Samuel showed that these functions agree with a polynomial $P_{I,M}(n)$ (called the \emph{Hilbert-Samuel polynomial}) of degree $d=\dim M$ for $n$ sufficiently large.  We can always write $P_{I,M}(n)$ in the form
\[P_{I,M}(n)=\sum_{i=0}^{d}(-1)^{i}e_{i}(I, M)\binom{n+d-i-1}{d-i}\] 
for unique numbers $e_{i}(I,M)$, known as the \emph{Hilbert coefficients} for $I$ with respect to $M$.  The largest number for which $H_{I,M}(n)$ and $P_{I,M}(n)$ disagree is called the  \emph{postulation number} for $I$, denoted $n(I,M):=\min\{j \  | \  H_{I,M}(n) = P_{I,M}(n) \  \forall n>j\}$.  Whenever $M=R$ we often suppress the $M$.

In this note, we focus our attention on non-positivity of the Hilbert coefficients for a parameter ideal.  Recall that a ring $R$ is \emph{unmixed} if $\dim \hat{R}/p=\dim R$ for all $p\in \Ass_{\hat{R}}\hat{R}$ where $\hat{R}$ denotes the $\m$-adic completion of $R$.  Our work was inspired by the following result of Ghezzi, et al \cite{GGHOPV} which characterizes the Cohen-Macaulayness of a ring in terms of the first Hilbert coefficient of a parameter ideal.  By a parameter ideal, we mean an ideal generated by a full system of parameters.

\begin{thm} \cite{GGHOPV}
Suppose $(R,\m)$ is an unmixed local ring and $q$ a parameter ideal. Then $e_{1} (q)\leq 0$
with equality if and only if R is Cohen-Macaulay.
\end{thm}

With the assumption that $\depth R \geq d-1$, we are able to prove the following:
\begin{thm}
Let $(R, \m)$ be a local Noetherian ring of dimension $d\geq 2$.  Suppose that $\depth R \geq d-1$.  If $q$ is a parameter ideal of $R$, then the following hold:
\begin{enumerate}
\item $e_{2}(q) \leq 0$
\item$ e_{2}(q)=0$  if and only if $n(q) < 2-d$ and $\grade gr_{q}(R)_{+} \geq d-1$
\item $e_{2}(q)=0$ implies $e_{3}(q)=e_{4}(q)= \cdots =e_{d}(q)=0$.
\end{enumerate}
\end{thm}

Here, $gr_{q}(R)$ denotes the associated graded ring of $R$ with respect to $q$.

We also discuss some results with respect to the first difference function, $\Delta$, defined in Section \ref{sec: other coeffs} and use this to examine the other Hilbert coefficients of a parameter ideal under the additional assumption that $\depth gr_{q}(R) \geq d-1$.

\section{The Hilbert Coefficients in Dimension One}
In this section we will discuss some results on the Hilbert coefficients of a parameter ideal in a one dimensional ring.

\begin{defn}
\textrm{Let $f:\mathbb{Z} \ra \mathbb{Z}$.  The first difference function, $\Delta(f)$, is defined by $\Delta(f(n))=f(n+1)-f(n)$.  We define the $i^{th}$ difference function inductively by $\Delta^{i}(f)=\Delta(\Delta^{i-1}(f))$.}  By convention, we define $\Delta^{0}(f)=f$.
\end{defn}

We first give a formula for the Hilbert coefficients in a one-dimensional ring.
\begin{prop} \label{lem: dim 1 hilbert fnct}
Suppose $(R, \m)$ is a one-dimensional local Noetherian ring and $q=(x)\subseteq R$ is a parameter ideal.  Then $((x^{i+1}):x^{i})=((x^{i+2}):x^{i+1})$ for all $i \gg 0$.  We set $l=\min\{i \ | \ ((x^{n+1}):x^{n})=((x^{i+1}):x^{i}) \textrm{ for all }n \geq i\}$ and $\tilde{x}=((x^{l+1}):x^{l})$. Then
\begin{enumerate}
\item \begin{enumerate} \item $e_{0}(q)=\lR(R/\tilde{x})$, and
\item $e_{1}(q)=\sum_{i=0}^{l-1}\left(\lR(R/\tilde{x})-\lR(R/((x^{i+1}):x^{i}))\right)$ for a fixed integer $l$.

\end{enumerate}

\item \begin{enumerate} \item $P_{q}(n)-H_{q}(n)=\sum_{i=n}^{\infty}\left(\lR(R/((x^{i+1}):x^{i}))-\lambda(R/\tilde{x})\right)$, and
\item $P_{q}(n) \geq H_{q}(n)$ for all $n\geq 0$.
\end{enumerate}
\end{enumerate}

\begin{proof}
Write $q=(x)$.  Note that $\lR((x^{i})/(x^{i+1}))=\lR(R/((x^{i+1}):x^{i}))$ for all $i$ as \newline $((x^{i+1}):x^{i})$ is the kernel of the surjective map $R \rightarrow (x^{i})/(x^{i+1})$ defined by $1 \mapsto \bar{x^{i}}$.  Then $\lambda(R/q^{n})=\sum_{i=0}^{n-1}\lR((x^{i})/(x^{i+1}))=\sum_{i=0}^{n-1}\lR(R/((x^{i+1}):x^{i}))$.  Note the ascending chain $$((x):x^{0}) \subseteq  ((x^{2}):x) \subseteq ((x^{3}):x^{2}) \subseteq \cdots$$ must stabilize.  Let $$l=\min\{i \ | \ ((x^{n+1}):x^{n})=((x^{i+1}):x^{i}) \textrm{ for all }n \geq i\}$$ and set $\tilde{x}=((x^{l+1}):x^{l})$.

For $n \geq l$, we have $\lR(R/(x^{n}))=\sum_{i=0}^{l-1}\lR(R/((x^{i+1}):x^{i}))+(n-l)\lR(R/\tilde{x})$.  This gives that \[P_{q}(n)=\sum_{i=0}^{l-1}\lR(R/((x^{i+1}):x^{i}))+(n-l)\lR(R/\tilde{x}).\]
From this, we see that $e_{0}(q)=\lR(R/\tilde{x})$ and $e_{1}(q)=\sum_{i=0}^{l-1}\left[\lR(R/\tilde{x})-\lR(R/((x^{i+1}:x^{i}))\right]$.  This proves (1).

Now if $n\leq l-1$, then $H_{q}(n)=\sum_{i=0}^{n-1}\lR(R/((x^{i+1}):x^{i}))$, and
\begin{eqnarray*}
P_{q}(n)-H_{q}(n)&=&\sum_{i=n}^{l-1}\lR(R/((x^{i+1}):x^{i}))+(n-l)\lR(R/\tilde{x})\\
&=&\sum_{i=n}^{l-1}\left(\lR(R/((x^{i+1}):x^{i}))-\lR(R/\tilde{x})\right) \\
&=& \sum_{i=n}^{\infty}\left(\lR(R/((x^{i+1}):x^{i}))-\lR(R/\tilde{x})\right),
\end{eqnarray*}
where the last equality holds since $((x^{i+1}):x^{i})=\tilde{x}$ for all $i \geq l$.  This gives 2($a$).

Note that for all $i$, we have $((x^{i+1}):x^{i}) \subset \tilde{x}$, so $\lR(R/((x^{i+1}):x^{i})) \geq \lR(R/\tilde{x})$ and we have $P_{q}(n)-H_{q}(n) \geq 0$.  In fact, if $n\geq l$, we have $P_{q}(n)-H_{q}(n)=0$.  This gives part 2($b$) of the proposition.
\end{proof}
\end{prop}

This proposition gives us a formula for the postulation number, $n(q)$,  of a parameter ideal $q$ in a one-dimensional ring.

\begin{cor}
Let  $(R, \m)$ be a one-dimensional local Noetherian ring and $q=(x)$ a parameter ideal.  Then $$n(q)=\min\{i \ | \ ((x^{i+1}):x^{i})=((x^{j+1}):x^{j}) \text{ for all } j\geq i \}-1.$$
\begin{proof}
Let $l=\min\{i \ | \ ((x^{i+1}):x^{i})=((x^{j+1}):x^{j}) \text{ for all } j \geq i \}$ and $\tilde{x}=((x^{l+1}):x^{l})$.  Then, using part 2($a$) of Proposition \ref{lem: dim 1 hilbert fnct}, clearly $n(q) \leq l-1$.  If $n(q) < l-1$, then we have $P_{q}(l)=H_{q}(l)$ and using 2($a$) again, this gives $\lR(R/((x^{l}):x^{l-1}))=\lR(R/\tilde{x})$.  But this contradicts the minimality of $l$.  Thus, we must have $n(q)=l-1$.
\end{proof}
\end{cor}

\begin{cor} \label{cor: dim 1, e_{1}}
Let $(R, \m)$ be a one-dimensional local Noetherian ring and $q=(x)$ a parameter ideal.  Then
\begin{enumerate}
\item For $k \in \mathbb{Z}$, if $P_{q}(k)-H_{q}(k)=0$, then $P_{q}(n)-H_{q}(n)=0$ for all $n \geq k$, i.e., $k>n(q)$.
\item $\Delta^{2}(P_{q}(n)-H_{q}(n))=\lR(((x^{n+2}):x^{n+1})/((x^{n+1}):x^{n}))$ for all $n$.
\end{enumerate}
\begin{proof}

For the first statement, suppose $P_{q}(k)-H_{q}(k)=0$.  Let $\tilde{x}$ be defined as in Proposition \ref{lem: dim 1 hilbert fnct}.  Then $P_{q}(k)-H_{q}(k)=\sum_{i=k}^{\infty}\left(\lR(R/((x^{i+1}):x^{i})-\lR(R/\tilde{x})\right)=0$.  Since \newline $\lR(R/((x^{i+1}):x^{i}))-\lR(R/\tilde{x})\geq 0$ for all $i\geq 0$, we must have equality for each $i\geq k$.  It follows that $P_{q}(n)-H_{q}(n)=0$ for all $n \geq k$, i.e., $k > n(q)$.

For (2), note by Proposition \ref{lem: dim 1 hilbert fnct}
\begin{eqnarray*}
\Delta(P_{q}(n)-H_{q}(n))&=&\lR(R/\tilde{x})-\lR(R/((x^{n+1}):x^{n})).
\end{eqnarray*}
So,
\begin{eqnarray*}
\Delta^{2}(P_{q}(n)-H_{q}(n))&=&\Delta(\Delta(P_{q}(n)-H_{q}(n))) \\
&=&\Delta(\lR(R/\tilde{x})-\lR(R/((x^{n+1}):x^{n}))) \\
&=&-\lR(R/((x^{n+2}):x^{n+1}))+\lR(R/((x^{n+1}):x^{n}))\\
&=& \lR(((x^{n+2}):x^{n+1})/((x^{n+1}):x^{n})). 
\end{eqnarray*}
In particular, this shows $\Delta(P_{q}(n)-H_{q}(n)) \leq 0$ and $\Delta^{2}(P_{q}(n)-H_{q}(n))\geq 0$ for all $n$.
\end{proof}
\end{cor}

\section{ The Second Hilbert Coefficient}

When working with Hilbert functions, a common technique is to reduce by a superficial sequence to obtain a ring of smaller dimension.  The following Proposition due to Nagata guarantees that when we do this, the Hilbert coefficients behave nicely.

\begin{prop} [cf. \cite{NA}, 22.6] \label{NA}
Let $(R, \m)$ be a Noetherian local ring, $I$ an $\m$-primary ideal, and $M$ a nonzero finitely generated $R$-module of dimension $d$.  Suppose $y\in $I is superficial with respect to $M$.  Then $\lR(0:_{M}y)$ is finite and 
\[P_{\bar{I}, \bar{M}}(n)=P_{I,M}(n)-P_{I,M}(n-1)+\lR(0:_{M}y).\]
In particular, we have
\[ e_{i}(\bar{I}, \bar{M})= \begin{cases}
e_{i}(I, M) & \textrm{for } i=0, \dots, d-2 \\
e_{d-1}(I, M)+(-1)^{d-1}\lR(0:_{M}y) & \textrm{for } i=d-1. \\
\end{cases} \]
\end{prop}

If $x \in I \backslash I^{2}$ is a non-zero-divisor and $x^{*}$ is a regular element of $gr_{I}(R)$, it can be easily shown that $n(I/(x))=n(I)+1$, so the postulation number also behaves nicely when we reduce via a superficial non-zero-divisor.

We begin with a formula for the last Hilbert coefficient.
\begin{lemma} \label{lem: e_{d} formula}
Suppose $(R,\m)$ has dimension $d$.  Let $I$ be an $\m$-primary ideal and $y\in I$ a superficial element.  Let $\bar{I}=I/(y)$, $H_{\bar{I}}(k)=\lR(R/(I^{k}, y))$ and $P_{\bar{I}}(k)$ denote the Hilbert Samuel polynomial for $\bar{I}$.  Then for $l\gg0$,
\[(-1)^{d}e_{d}(I)=\sum_{k=1}^{l}\left(H_{\bar{I}}(k)-P_{\bar{I}}(k)\right) -
\sum_{k=1}^{l} \lR((I^{k}:y)/I^{k-1})+l\lR(0:y).\]
Furthermore, if $y$ is also a non-zero-divisor on $R$, we have
\[(-1)^{d}e_{d}(I)=\sum_{k=1}^{\infty}\left(H_{\bar{I}}(k)-P_{\bar{I}}(k)\right) -
\sum_{k=1}^{\infty} \lR((I^{k}:y)/I^{k-1}).\]
\begin{proof}
For $k \in \mathbb{Z}$, consider the exact sequence:
\[0 \rightarrow \frac{I^{k}:y}{I^{k-1}} \rightarrow R/I^{k-1} \xrightarrow{y} R/I^{k}
\rightarrow R/(I^{k},y) \rightarrow 0. \]
From this we see that
$\lR(R/(y,I^{k}))=\lR(R/I^{k})-\lR(R/I^{k-1})+\lR\left(\frac{I^{k}:y}{I^{k-1}}\right)$.
Subtracting $P_{\bar{I}}(k)$ and summing both sides, we get, for
$l\gg0$,
\begin{eqnarray*}
\sum_{k=1}^{l}\left(\lambda(R/(y,
I^{k})-P_{\bar{I}}(k)\right)&=&\sum_{k=1}^{l}
\left(\lambda(R/I^{k})-\lambda(R/I^{k-1})+\lambda\left(\frac{I^{k}:y}{I^{k-1}}\right)-
P_{\bar{I}}(k)\right). 
\end{eqnarray*}
Thus, \begin{eqnarray*}
\sum_{k=1}^{l}(H_{\bar{I}}(k)-P_{\bar{I}}(k))&=&
\lambda(R/I^{l})-\sum_{k=1}^{l}\sum_{i=0}^{d-1}(-1)^{i}
\binom{k+d-2-i}{d-1-i}e_{i}(\bar{I})+\sum_{k=1}^{l}\lambda \left(\frac{I^{k}:y}{I^{k-1}}\right)\\
&=& \sum_{i=0}^{d}(-1)^{i}\binom{l+d-1-i}{d-i}
e_{i}(I)-\sum_{i=0}^{d-1}(-1)^{i}
\binom{l+d-1-i}{d-i}e_{i}(\bar{I})\\
&& \hspace{.5cm}+\sum_{k=1}^{l}\lambda\left(\frac{I^{k}:y}{I^{k-1}}\right)
\end{eqnarray*}
where $\lambda(-)=\lR(-)$.

By
Proposition \ref{NA}, we have $e_{i}(I)=e_{i}(\bar{I})$
for $i=0, \dots , d-2$ and
$e_{d-1}(I)=e_{d-1}(\bar{I})-(-1)^{d-1}\lR(0:y)$.  Hence,
\begin{eqnarray*}
\sum_{k=1}^{l}(H_{\bar{I}}(k)-P_{\bar{I}}(k))&=&-l\lR(0:y)+(-1)^{d}e_{d}(I)+\sum_{k=1}^{l}\lR((I^{k}:y)/I^{k-1}). 
\end{eqnarray*}
Rearranging, we get
\[(-1)^{d}e_{d}(I)=\sum_{k=1}^{l}\left(H_{\bar{I}}(k)-P_{\bar{I}}(k)\right) - \sum_{k=1}^{l} \lR((I^{k}:y)/I^{k-1})+l\lR(0:y)\]
and if $y$ is also a non-zero-divisor on $R$, we have
\[(-1)^{d}e_{d}(I)=\sum_{k=1}^{\infty}\left(H_{\bar{I}}(k)-P_{\bar{I}}(k)\right) - \sum_{k=1}^{\infty} \lR((I^{k}:y)/I^{k-1})\]
since for $k\gg0$, $H_{\bar{I}}(k)-P_{\bar{I}}(k)=0$ and $\lR((I^{k}:y)/I^{k-1})=0$.
\end{proof}
\end{lemma}

We define $\grade gr_{I}(R)_{+}$ to be the maximal length of a regular sequence for $gr_{I}(R)$ contained in $gr_{I}(R)_{+}$.  Then $\grade gr_{I}(R)_{+}=\depth gr_{I}(R)$.  For $x\in I^{n}\backslash I^{n+1}$ let $x^{*}$ denote its image in $I^{n}/I^{n+1} \subseteq gr_{I}(R)$.  The grade of the associated graded ring also behaves nicely with
respect to superficial sequences as evidenced by the following
lemmas.  Lemma \ref{lem: HM 2.2} is also known as ``Sally's Machine".

\begin{lemma} \cite[Lemma 2.1]{HM} \label{lem: HM 2.1} Let $x_{1}, \dots, x_{k}$ be a superficial
sequence for $I$.  If $\grade gr_{I}(R)_{+} \geq k$, then $x_{1}^{*},  \dots , x_{k}^{*}$ is a regular
sequence.
\end{lemma}

\begin{lemma}\cite[Lemma 2.2]{HM} \label{lem: HM 2.2}
Suppose $y_{1}, \dots, y_{k}$ is a superficial sequence for an ideal $I$.  Let $\bar{R}$ and
$\bar{I}$ denote $R/(y_{1}, \dots, y_{k})$ and $I/(y_{1}, \dots, y_{k})$, respectively. If
$\grade gr_{\bar{I}}(\bar{R})_{+} \geq 1$, then $\grade gr_{I}(R)_{+} \geq k+1$.
\end{lemma}

We now consider the second Hilbert coefficient, $e_{2}(q)$ for a parameter ideal $q$.

\begin{thm} \label{2nd coeff}
Suppose $(R, \m)$ is a Noetherian local ring of dimension $d \geq 2$ and $\depth R \geq d-1$.  Let $q \subseteq R$ be a parameter ideal.  Then
\begin{enumerate}
\item $e_{2}(q) \leq 0.$
\item $e_{2}(q)=0$ if and only if $n(q)<2-d$ and $\depth gr_{q}(R) \geq d-1.$
\item $e_{2}(q)=0$ implies $e_{3}(q)=\cdots = e_{d}(q)=0$.
\end{enumerate}
\begin{proof}
We may assume that $R$ has infinite residue field by passing to $R[x]_{\m R[x]}$ if necessary.  We will proceed by induction on $d=\dim R$.  First suppose $d=2$.  Let $q=(y,x)$ where $y \in q\backslash \m q$ is a superficial non-zero-divisor for $R$.  Let $(\overline{\cdot })$ denote working modulo $(y)$.  Now, $\bar{q}$ is a parameter ideal in the one-dimensional ring $\bar{R}$, so by Proposition \ref{lem: dim 1 hilbert fnct}, $H_{\bar{q}}(k)-P_{\bar{q}}(k) \leq 0$ for all $k \geq 0$.  In particular, Lemma \ref{lem: e_{d} formula} gives
\begin{eqnarray*}
e_{2}(q) & = & \sum_{k=1}^{\infty}(H_{\bar{q}}(k)-P_{\bar{q}}(k)) - \sum_{k=1}^{\infty} \lR((q^{k}:y)/q^{k-1})\\
& \leq & 0.
\end{eqnarray*}
Note that if the left-hand side of the equation above is zero, we must have that $\lR((q^{k}:y)/q^{k-1})=0$ and $P_{\bar{q}}(k)=H_{\bar{q}}(k)$ for all $k\geq 1$.  In particular, the condition $\lR((q^{k}:y)/q^{k-1})=0$ for all $k \geq 1$ implies that $y^{*}$ is a non-zero-divisor in $gr_{q}(R)$, so $\depth gr_{q}(R) \geq 1$.  Now, since $y^{*}$ is a non-zero-divisor, $n(\bar{q})=n(q)+1$; i.e., $n(q)<0$.  This proves the Corollary when $d=2$.

Now if $\dim R >2$, then let $y_{1}, \dots, y_{d-2} \in q\backslash \m q$ be a superficial sequence of non-zero-divisors for $R$.  Then $\bar{q}=q/(y_{1}, \dots, y_{d-2})$ is a parameter ideal in the two-dimensional ring $\bar{R}=R/(y_{1}, \dots, y_{d-2})$ which has $\depth \bar{R} \geq 1$.  Hence, by induction, we have $e_{2}(q)=e_{2}(\bar{q}) \leq 0$.

For (2), first suppose $e_{2}(q)=0$.  Then by induction $\grade gr_{\bar{q}}(\bar{R})_{+} \geq 1$.  By Lemma \ref{lem: HM
2.2}, this implies $\grade gr_{q}(R)_{+} \geq d-2+1=d-1.$ Finally,
this gives $y_{1}^{*}, \dots, y_{d-2}^{*}$ is a regular sequence by Lemma \ref{lem: HM 2.1}. Hence,
$n(\bar{q})< 0$ if and only if $n(q) < 2-d$. This
gives the forward implications for (2).

For the backward implication of (2), suppose $n(q)<2-d$ (i.e., $H_{q}(n)=P_{q}(n)$ for all $n \geq 2-d$) and $\grade gr_{q}(R)_{+} \geq d-1$.  Then $P_{q}(n)=0$ for all $2-d \leq n \leq 0$.  Plugging the values $0, -1, -2, \dots, 2-d$ successively
into $P_{q}(n)$, one can see that we get $e_{d}(q)=e_{d-1}(q)=\cdots=e_{2}(q)=0$. 

Finally, (3) follows from the proof of (2).
\end{proof}
\end{thm}

\begin{cor}
Suppose $(R, \m)$ is a local Noetherian ring of dimension $d\geq 2$ and $\depth R \geq d-1$.  Then for any parameter ideal $q \subseteq R$, we have $$\lR(R/q) \leq e_{0}(q)-e_{1}(q).$$
\begin{proof} As before, we may assume that $R/\m$ is infinite.
From Proposition \ref{lem: dim 1 hilbert fnct}, we have that $H_{\bar{q}}(n) \leq P_{\bar{q}}(n)$ for all $n \geq 1$, where $\bar{q}=q/(y_{1}, \dots, y_{d-1})$ for $y_{1}, \dots, y_{d-1} \in q$ a superficial sequence which is part of a minimal generating set for $q$.  Note that we may also choose $y_{1}, \dots, y_{d-1}$ to be a regular sequence as $\depth R \geq d-1$.  Now, letting $n=1$ and using the fact that $e_{i}(q)=e_{i}(\bar{q})$ for $i=0,1$ since $y_{1}, \dots, y_{d-1}$ is a superficial and regular sequence, the result follows.
\end{proof}
\end{cor}

The assumption that $\depth R \geq d-1$ is necessary in Theorem \ref{2nd coeff}, as evidenced by the following example.  We use Macaulay2 \cite{M2} to compute the example.

\begin{example}
Let $R=k[x,y,z,u,v,w]/I$ where $I$ is the intersection of ideals $I=(x+y, z-u, w)\cap (z, u-v, y) \cap (x,u,w)$ and $q=(u-y, z+w, x-v)$.  Then $R$ is an unmixed ring of dimension three and depth one and $q$ is a parameter ideal with \[P_{q}(n)=3\binom{n+2}{3}+2\binom{n+1}{2}+n.\]  In particular, $e_{2}(q)=1> 0.$
\end{example}

Note that in the example above, one could mod out the ring $R$ by a superficial non-zero-divisor in $q \backslash \m q$ to obtain an example of a two-dimensional ring $\bar{R}$ of depth zero with parameter ideal $\bar{q}$ satisfying $e_{2}(\bar{q})=1>0$.

The upper bound for $e_{2}(q)$ in Theorem \ref{2nd coeff} can be achieved even if $R$ is not Cohen-Macaulay.  We also provide an example below with negative second Hilbert coefficient.  In both examples, we use the software system Macaulay2 \cite{M2} to compute the Hilbert-Samuel functions.

\begin{example}
Let $R=k[[x^{5}, xy^{4}, x^{4}y, y^{5}]] \cong k[[t_{1}, t_{2}, t_{3}, t_{4}]]/J$ where $J$ is the ideal $J=(t_{2}t_{3}-t_{1}t_{4}, t_{2}^{4}-t_{3}t_{4}^{3}, t_{1}t_{2}^{3}-t_{3}^{2}t_{4}^{2}, t_{1}^{2}t_{2}^{2}-t_{3}^{3}t_{4}, t_{1}^{3}t_{2}-t_{3}^{4}, t_{3}^{5}-t_{1}^{4}t_{4})$.  Then $R$ is a two-dimensional complete domain with depth one.  The parameter ideal $q=(x^{5}, y^{5})$ has Hilbert-Samuel polynomial $$P_{q}(n)=5\binom{n+1}{2}+2n,$$ so $e_{2}(q)=0$.
\end{example}

\begin{example}
Let $R=k[x,y,z,t]/((x^{2}, z^{4}) \cap (x-y, z+t))$.  Then $R$ is a two-dimensional unmixed ring with depth one.   The ideal $q=(x+t+y, z-y)$ is a parameter ideal with Hilbert-Samuel polynomial $$P_{q}(n)=9 \binom{n+1}{2} +2 \binom{n}{1}-1.$$  Hence, $e_{2}(q)=-1<0$.
In this example, we have that $n(q)=0$; that is, $P_{q}(0)\neq H_{q}(0)$ and $P_{q}(n)=H_{q}(n)$ for all $n \geq 1$.  However, we do have that $\depth gr_{q}(R)\geq 1$.  \end{example}

\section{The Higher Hilbert Coefficients} \label{sec: other coeffs}
In our first theorem of this section we use techniques similar to those of Marley to obtain a result reminiscent of Theorem 1 in \cite{MA}.
\begin{thm} \label{thm: delta fnct}
Let $(R, \m)$ be a local Noetherian ring of dimension $d$ and suppose $q$ is a parameter ideal for $R$ satisfying $\depth gr_{q}(R)\geq d-1$.  Then for $0\leq i \leq d+1$ and $n \in \mathbb{Z}$,  \[(-1)^{i}\Delta^{d+1-i}(P_{q}(n)-H_{q}(n)) \geq 0.\]
\begin{proof}
We first note that it is enough to prove the result when $i=0$.  Indeed, suppose $g:\mathbb{Z}\ra\mathbb{Z}$ satisfies $g(n)=0$ for all $n$ sufficiently large and $\Delta(g(n))\geq 0$ for all $n.$  
  Then we claim $g(n)\leq 0$ for all $n$.  Let $N$ be such that $g(n)=0$ for all $n \geq N$. Then $\Delta(g(N-1))=g(N)-g(N-1)\geq 0$ implies $g(N-1) \leq 0$.  Inductively, one can show that $g(j) \leq 0$ for all $j$.  In particular, if we set $g(n)=P_{q}(n)-H_{q}(n)$, and assume $(-1)^{i}\Delta^{d+1-i}(g(n)) \geq 0$ for all 
$n$, then $(-1)^{i}\Delta^{d+1-i-1}(g(n)) \leq 0$ gives the theorem for $i+1$.  Hence, it is enough to prove \[\Delta^{d+1}(P_{q}(n)-H_{q}(n)) \geq 0 \ \ \ \text{for all }n.\]

We will use induction on the dimension $d$.  Note the case $d=1$ is proved in Corollary \ref{cor: dim 1, e_{1}}.  Suppose $d>1$.  Let $a\in q\backslash \m q$ such that $a^{*}$ is a $gr_{q}(R)$-regular element.  Let $\bar{q}=q/(a)$ and $\bar{R}=R/(a)$.  Then note that $\depth _{\bar{q}}(\bar{R}) \geq d-2$ and $\bar{q}$ is a parameter ideal for the $(d-1)$-dimensional ring $\bar{R}$.  So, by induction, \[\Delta^{d}(P_{\bar{q}}(n)-H_{\bar{q}}(n))\geq 0 \ \ \ \text{for all } n.\]
Now, as $a^{*}$ is a non-zero-divisor in $gr_{q}(R)$, we have $H_{\bar{q}}(n)=H_{q}(n)-H_{q}(n-1)$ for all $n$.  Similarly, $P_{\bar{q}}(n)=P_{q}(n)-P_{q}(n-1)$.  Hence,
\begin{eqnarray*}
\Delta^{d+1}(P_{q}(n)-H_{q}(n))&=&\Delta^{d}(\Delta(P_{q}(n)-H_{q}(n)) \\
&=&\Delta^{d}(P_{\bar{q}}(n+1)-H_{\bar{q}}(n+1)) \\
&\geq&0 \ \ \ \ \text{for all } n. 
\end{eqnarray*}
Thus, \[\Delta^{d+1}(P_{q}(n)-H_{q}(n))\geq0 \ \ \ \ \text{for all } n. \qedhere\]
\end{proof}

\end{thm}

\begin{cor} \label{cor: difference stabilizes}
Let $(R, \m)$ be a local Noetherian ring of dimension $d$ and suppose $q$ is a parameter ideal for $R$ satisfying $\depth gr_{q}(R)\geq d-1$.  Suppose $P_{q}(k)-H_{q}(k)=0$ for some $k$.  Then $P_{q}(n)-H_{q}(n)=0$ for all $n\geq k$, i.e., $k>n(q)$.
\begin{proof}
Letting $i=d$ in Theorem \ref{thm: delta fnct}, we have $(-1)^{d}\Delta(P_{q}(n)-H_{q}(n)) \geq 0$ for all $n$.  This gives $(-1)^{d}(P_{q}(n+1)-H_{q}(n+1)) \geq (-1)^{d}( P_{q}(n)-H_{q}(n))$ for all $n$.  In particular, we have \[0=(-1)^{d}( P_{q}(k)-H_{q}(k))\leq(-1)^{d}(P_{q}(n)-H_{q}(n))\leq 0 \ \ \ \forall \ n \geq k \]

where the last inequality holds because $P_{q}(N)-H_{q}(N)=0$ for $N\gg0$.  Thus, $P_{q}(n)=H_{q}(n)$ for all $n \geq k$.
\end{proof}

\end{cor}

\begin{remark} \label{rmk: post and zero}
Let $(R, \m)$ be a local Noetherian ring of dimension $d$ and suppose $q$ is a parameter ideal for $R$.   For $0\leq i \leq d-1$, we have the following:
\begin{enumerate}
\item If $n(q)<i-d$, then $e_{j}(q)=0$ for $j \geq i$.
\item If $\depth gr_{q}(R) \geq d-1$, the converse to (1) holds.
\end{enumerate}
\begin{proof}
Note (1) follows by using the fact that $P_{q}(j)=0$ for $i-d<j<0$.  For (2), suppose  $\depth gr_{q}(R) \geq d-1$ and $e_{i}(q)=0$ for $j \geq i$.  Then $P_{q}(i-d)=0=H_{q}(i-d)$ and by Corollary \ref{cor: difference stabilizes}, $n(q)<i-d$.
\end{proof}
\end{remark}

\begin{ques}
Does the converse to part (1) of Remark \ref{rmk: post and zero} above hold in general?
\end{ques}

\begin{cor} \label{cor: e_{i} nonpositive}
Let $(R, \m)$ be a local Noetherian ring of dimension $d$ and suppose $q$ is a parameter ideal for $R$ satisfying $\depth gr_{q}(R)\geq d-1$.   Then for $1\leq i \leq d$
\begin{enumerate}
\item $e_{i}(q) \leq 0.$
\item $(-1)^{j+1}(e_{0}(q)-e_{1}(q)+\cdots+(-1)^{j}e_{j}(q)-\lR(R/q))\geq 0$ for $j=1, \dots, d.$
\end{enumerate}
\begin{proof}
Note that it is enough to prove (1) in the case $i=d$, as we can then use reduction by a superficial sequence to obtain $e_{i}(q)\leq0$ for $i=1, \dots, d-1$.  Letting $i=d+1$ in Theorem \ref{thm: delta fnct}, we have
\begin{equation}
(-1)^{d+1}(P_{q}(n)-H_{q}(n))\geq 0\ \ \ \text{for all } n. \label{delta i=d+1}
\end{equation}
  If $n=0$, $(-1)^{d+1}((-1)^{d}e_{d}(q)-H_{q}(0))\geq 0$ implies $-e_{d}(q)\geq 0$; that is, $e_{d}(q)\leq 0$.

For (2), we will first prove the case $j=\dim R=d$.  Indeed, letting $n=1$ in equation $\text{(\ref{delta i=d+1})}$, we see \[(-1)^{d+1}(e_{0}(q)-e_{1}(q)+\cdots +(-1)^{d}e_{d}-\lR(R/q))\geq 0.\]
Now, let $a_{1}, \dots a_{d-j} \in q\backslash q^{2}$ be part of a minimal generating set for $q$ such that $a_{1}^{*}, \dots, a_{d-j}^{*}$ is a $gr_{q}(R)$-regular sequence.  Then, setting $\bar{R}=R/(a_{1}, \dots a_{d-j} )$ and $\bar{q}=q/(a_{1}, \dots a_{d-j} )$ we have $\bar{q}$ is a parameter ideal in the $j$-dimensional ring $\bar{R}$, and $\depth gr_{\bar{q}}(\bar{R})\geq j-1$.  Finally, $\lR(R/q)=\lambda_{\bar{R}}(\bar{R}/\bar{q})$ and since $a_{1}, \dots a_{d-j} $ defines a superficial regular sequence in $R$, we have $e_{i}(\bar{q})=e_{i}(q)$ for all $i=0, \dots, j$.  It follows that  \[(-1)^{j+1}(e_{0}(q)-e_{1}(q)+\cdots+(-1)^{j}e_{j}(q)-\lR(R/q))\geq 0. \qedhere \]
\end{proof} 
\end{cor}

\begin{cor} \label{cor: e_{j}=0}
Let $(R, \m)$ be a local Noetherian ring of dimension $d$ and suppose $q$ is a parameter ideal for $R$ satisfying $\depth gr_{q}(R)\geq d-1$.  Suppose $e_{i}(q)=0$ for some $1\leq i\leq d-1$.  Then $e_{j}(q)=0$ for $i\leq j\leq d$.
\begin{proof}
Note that it is enough to prove that $e_{i+1}(q)=0$.  Reducing by a superficial sequence if necessary, we may assume that $i=d-1$.  Since $e_{0}(q)>0$, we must have that $d>1$, so by assumption, $\depth gr_{q}(R) >0$.  Let $a\in q$ be such that $a^{*}\in gr_{q}(R)$ is a non-zero-divisor.  Then $e_{d-1}(\bar{q})=e_{d-1}(q)=0$ implies that $P_{\bar{q}}(0)=0=H_{\bar{q}}(0)$.  Now, by Corollary \ref{cor: difference stabilizes}, $n(\bar{q})\leq -1$.  As $n(\bar{q})=n(q)+1$, this gives $n(q)\leq -2$, and in particular, $(-1)^{d}e_{d}(q)=P_{q}(0)=H_{q}(0)=0$.
\end{proof}
\end{cor}

\section*{Acknowledgments}
The work in this paper was done in preparation for the author's dissertation and would not have been possible without the guidance, encouragement, and input from my advisor Tom Marley.  The author is also grateful to Brian Harbourne and Luchezar Avramov for their feedback on this research.


\begin{thebibliography}{9}
\bibitem{BH} W. Bruns and J. Herzog, \textit{Cohen-Macaulay Rings}, Cambridge University Press, Cambridge, 1993.

\bibitem{GGHOPV} L. Ghezzi, S. Goto, J. Hong, K. Ozeki, T.T. Phuong, and W.V. Vasconcelos. ``Negativity conjecture for the first Hilbert coefficient,'' preprint.

\bibitem{GHV} L. Ghezzi, J. Hong, and W.V. Vasconcelos, ``The signature of the Chern coefficients of local rings,''  \textit{Math. Research Letters}, \textbf{16} (2009), 279-289.



\bibitem{HM}  S. Huckaba and T. Marley, ``Hilbert coefficients and the depths of associated graded rings,'' \textit{J. London Math. Soc.} (2) 56 (1997), no. 1, 64-76.

\bibitem{HS} C. Huneke and I. Swanson, \textit{Integral Closure of Ideals, Rings, and Modules}, London Mathematical Society Lecture Note Series, 336. Cambridge University Press, Cambridge, 2006.

\bibitem{MA} T. Marley ``The coefficients of the Hilbert polynomial and the reduction number of an ideal,'' \textit{J. London Math. Soc.} (2) 40 (1989), no.1, 1-8.

\bibitem{M2} D. Grayson and M. Stillman, \textit{Macaulay2, a software system for research in algebraic geometry}, Available at http://www.math.uiuc.edu/Macaulay2.

\bibitem{MA1} H. Matsumura. \textit{Commutative Ring Theory}. Cambridge Studies in Advanced Mathematics, vol.~8, Cambridge University Press, Cambridge, 1986, Translated from the Japanese by M. Reid. \MR{MR879273 (88h:13001)}

\bibitem{NA} M. Nagata, \textit{Local Rings}, Interscience, New York, 1962.

\end{thebibliography}
\end{document}